\documentclass[twoside,12pt]{article}
\usepackage[all, 2cell, dvips]{xy}
\usepackage{latexsym,amsmath, amsfonts, graphics, epsf, epic}
\usepackage{amsthm}
\usepackage{amssymb}
\usepackage{amsopn}
\usepackage{amscd}
\usepackage{color}

\theoremstyle{plain}

\theoremstyle{definition}

\usepackage{amssymb}

\def\lm{\lambda}

\def\l.l.o.{\it l.l.o}

\def\chiup{\raise 2pt\hbox{$\chi$}}

\title{Bijections for an identity of SYT}

\begin{document}

\centerline{\bf Bijections for an identity of Young Tableaux}

\medskip
\centerline{A. Regev and D. Zeilberger}

\bigskip
Recall that  {\it partition} $\lambda=(\lambda_1, \dots, \lambda_k)$ is a weakly-decreasing sequence
of positive integers, and if $n=\lambda_1+ \dots + \lambda_k$ we say that $\lambda$ is a partition
of $n$ and write $\lm\vdash n$ and $l(\lambda)=k$.  For example $\lambda=(3,3,2,2)$ is a partition
of $10$ and $l(\lambda)=4$.

\medskip
A convenient way to represent a partition $\lambda$ is via its {\it Ferrers} diagram, that consists
of $l(\lambda)$ left-justified lines of dots, such that the $i$-th row has $\lambda_i$
dots. If you replace the dots by empty boxes, you would get what is called a
{\it Young} diagram (of shape $\lambda)$.

\medskip
Finally recall that a {\it standard Young tableau} (SYT)  of shape $\lambda \vdash n$
is any way of placing the integers $\{1,2, \dots, n\}$ into the empty boxes of
the Young diagram in such a way that all the rows and all the columns are increasing.
For example $[[1,2,4,6],[3,5,7],[8,9]]$ is a SYT of shape $(4,3,2)$.

\medskip
Let $H(k,\ell;n)=\{\lm\vdash n\mid \lm_{k+1}\le \ell\}$ denote the
partitions of $n$ in the $(k,\ell)$ hook. For example, $H(k,0;n)$
are the partitions of $\lm\vdash n$ with $\ell(\lm)\le k$. Let
$f^\lm$ denote the number of SYTs of
shape $\lm$. One then observes the following intriguing identity:
\[
\sum_{\mu\in H(1,1;n+1)} (f^\mu)^2=\sum_{\lm\in H(2,0;2n)}
f^\lambda.
\]
We give this identity a bijective proof by showing that both
\begin{eqnarray}\qquad\sum_{\mu\in H(1,1;n+1)} (f^\mu)^2,\end{eqnarray}
   \begin{eqnarray}\qquad\sum_{\lm\in H(2,0;2n)} f^\lambda,\end{eqnarray}
can be mapped {\it bijectively} to the set of row-increasing
matrices of shape $(n,n)$ whose set of entries is $\{1,2, \dots , 2n\}$, and hence, by composing, to each other.
We describe these bijections. We leave it as pleasant excercises to the reader to formally prove that
these are indeed bijections, by proposing inverse mappings,
and proving that the compositions (in both direction) yield the identity mapping in each case.

\medskip
 The input for both (1)
 and (2)  is
a $2\times n$ matrix of integers

$a_1 \ldots a_n$

$b_1 \ldots b_n$

\medskip
such that $\{a_1,\ldots,a_n\} \cup \{b_1, \ldots
b_n\}=\{1,2,\ldots, 2n\}$,

\medskip
$a_1<a_2<...<a_n$ and

$b_1<b_2<...<b_n$.

\bigskip
{\bf Description of the bijection for} (1)

\medskip
Here the output is: Two  standard tableaux of the same
$(1,1)$-hook shape. Let

$|\{a_1, \ldots, a_n\}\cap \{1, \ldots, n\}|=k,$ so $\{a_1,
\ldots, a_n\}\cap \{1, \ldots, n\}=\{a_{i_1}<\cdots < a_{i_k}\}.$
%
%
Form now a SYT in the $(1,1)$ hook as follows:

Its (first) row is \[~~1,~~a_{i_1}+1,~\cdots ~,a_{i_k}+1;\] its
(first) column is made of the remaining integers
\[
\{1,\ldots n+1\}\setminus \{~1,~a_{i_1}+1,~\cdots ~,a_{i_k}+1\}
=\{a'_{j_1},\ldots ,a'_{j_{n-k}}\} \]
 in increasing order. This gives the first SYT -- of $(1,1)$-hook shape
$(k+1,1^{n-k})$.

\bigskip


It follows that
\[
|\{a_1,\ldots, a_n\} \cap \{n+1, \ldots, 2n\}|= n-k,
\]
so denote
\[
\{a_1,\ldots, a_n\} \cap \{n+1, \ldots, 2n\}=
\{b_{t_1}<\cdots<b_{t_{n-k}}\}.
\]
Since $n+1\le b_{t_1}$, hence $2\le b_{t_1}-(n-1)$. Form the
numbers \[1<b_{t_1}-(n-1)<\cdots<b_{t_{n-k}}-(n-1)\] and place
them, in that (increasing) order, in the column of the second
(1,1)-hook shape tableau, and the complement integers
\[\{1,\ldots,n+1\}\setminus
\{1,b_{t_1}-(n-1),\ldots,b_{t_{n-k}}-(n-1)\}
\]
 in
increasing order, in the row (after the corner 1). This gives the
second SYT -- again of shape $(k+1,1^{n-k})$. This map is clearly
a bijection.

\medskip
{\bf Example of bijection (1).} Let $n=5$ and consider the
$2\times 5$ array
\begin{eqnarray*}
2,~4,~8,~9,~10\\ 1,~3,~5,~6,~7~
\end{eqnarray*}
Now $\{2,~4,~8,~9,~10\}\cap\{1,\ldots,5\}=\{2,4\}\to_{+1} \{3,5\}$
so the first row of the first tableau is $(1,3,5)$, hence
($n+1=6$) its first column is $(1,2,4,6)^T$.

\medskip
Also, $\{2,~4,~8,~9,~10\}\cap\{6,\ldots,10\}=\{8,9,10\}\to_{-4}
\{4,5,6\}$ so
 the first column of the second tableau is $(1,4,5,6)^T$, hence
 its first row is $(1,2,3)$. The pair of these two $(1,1)$-
 tableaux corresponds to the above array.

\bigskip
{\bf Description of the bijection for} (2)

\medskip
Here the output is a SYT whose shape is a $\leq 2$-rowed partition
of $2n$ If for all $i~ a_i<b_i$, then it is a SYT, and do nothing
(the output is the input).

\medskip
Otherwise, let $i$ be the smallest index such that $a_i>b_i.$
Replace the above array by

\newpage
$b_1 ~\ldots~b_{i-1} ~~b_i~~~\,a_i ~~~~a_{i+1}~~~\ldots~~~~a_n$

$a_1~\ldots \,~a_{i-1} ~b_{i+1}~b_{i+2}~b_{i+3}~\ldots~b_n$

\medskip

If this is a SYT, then stop. Otherwise continue: Typically we
arrive at a two rows array

\medskip
$c_1,~\ldots~c_s~\ldots~c_r$

$d_1~\ldots~~d_s$

\medskip
with $s\le r,~~r+s=2n$, with $~c_1<\cdots<c_r$ and with
$d_1<\cdots<d_s$. If $c_j<d_j,$ $~j=1,\ldots,s$ ~~then this array
is SYT and we are done. Otherwise let $i$ be the smallest index
such that $c_i>d_i$, then replace the above array by

\medskip
$d_1~\ldots~d_{i-1}~~d_i~~~\;c_i~~~~c_{i+1}~\ldots
c_{s-1}~~\ldots~~~c_r$

$c_1~\ldots~c_{i-1}~~d_{i+1}~\;d_{i+2}~d_{i+3}~\ldots~d_s$

\medskip
Continue until a SYT is reached: thus the process stops at a SYT
of shape $\lm\in H(2,0;2n)$.

\medskip
{\bf Example of bijection (2).} Again consider the same $2\times
5$ array as above, then the the bijection is as follows:
\begin{eqnarray*}
2,~4,~8,~9,~10~~~~\\ 1,~3,~5,~6,~7~\to
\end{eqnarray*}
\begin{eqnarray*}
1,~2,~4,~8,~9,~10~~\\ 3,~5,~6,~7~~~~\quad\to
\end{eqnarray*}
\begin{eqnarray*}
 3,~5,~6,~7,~8,~9,~10\\1,~2,~4~~~~~~~\qquad\to
\end{eqnarray*}
\begin{eqnarray*}
 1,~3,~5,~6,~7,~8,~9,~10\\2,~4~~~~~~~~~~~~~~~~~~~~~~~~~
\end{eqnarray*}

\medskip {\bf Remark.}
Let $\lm=(\lm_1,\lm_2,\ldots )$ be a partition,~~$\lm_1\ge
\lm_2\ge\ldots$ and denote

$\lm^{+1}=(\lm_1+1,\lm_2,\lm_3\ldots)$. The same bijections can be
applied to similar arrays of shape $(n+1,n)$, yielding a bijective
proof for the SYT identity
\[
\sum_{\mu\in H(2,0;2n+1)}f^\mu=\sum_{\lm\in H(1,1;n+2)}f^\lm\cdot
f^{\lm^{+1}}.
\]

\medskip {\bf Amitai Regev}, Department of Mathematics, The
Weizmann Institute of Science, Rehovot, Israel.
amitai.regev~at~weizmann.ac.il

\medskip

{\bf Doron Zeilberger}, Mathematics Department, Rutgers University
(New Brunswick), Piscataway, NJ, USA.
zeilberg~at~math.rutgers.edu .

\end{document}